\documentclass{amsart}
\usepackage{amsmath,array}
\usepackage[frame]{xy}
\usepackage{xspace}
\newtheorem{thm}{Theorem}[section]
\newtheorem{lem}[thm]{Lemma}
\newtheorem{cor}[thm]{Corollary}

\newtheorem{cnj}[thm]{Conjecture}

\theoremstyle{definition}
\newtheorem{dfn}[thm]{Definition}

\theoremstyle{remark}

\newtheorem*{eg*}{Example}

\newcommand{\Set}[1]{\ensuremath{\mathcal{#1}}}            
\newcommand{\Mat}[1]{\ensuremath{\mathbf{#1}}}             
\newcommand{\Dfn}[1]{\emph{#1}}                            

\let\size=\abs

\newcommand{\eval}[2][\right]{\relax
  \ifx#1\right\relax \left.\fi#2#1\rvert}

\newcommand{\shrink}[2][0pt]
{\hbox to #1{\hss #2\hss}}

\def\ddotsup{\mathinner{\mkern1mu\raise1pt
                        \vbox{\kern7pt\hbox{.}}
                        \mkern2mu\raise4pt\hbox{.}
                        \mkern2mu\raise7pt\hbox{.}\mkern1mu}}

\def\IsInteger#1{%
    TT\fi
    \begingroup \lccode`\-=`\0 \lccode`+=`\0
      \lccode`\1=`\0 \lccode`\2=`\0 \lccode`\3=`\0
      \lccode`\4=`\0 \lccode`\5=`\0 \lccode`\6=`\0
      \lccode`\7=`\0 \lccode`\8=`\0 \lccode`\9=`\0
    \lowercase{\endgroup
      \expandafter\ifx\expandafter\delimiter
      \romannumeral0\string#1}\delimiter
  }

\newcounter{seq}
\newcommand{\seq}[4][,]
{\if\IsInteger{#3}%
{\setcounter{seq}{#3}%
\ifmmode%
#2{\arabic{seq}}#1\stepcounter{seq}#2{\arabic{seq}}#1\dots #1#2{#4}%
\else%
$#2{\arabic{seq}}$#1~\stepcounter{seq} $#2{\arabic{seq}}$#1 \dots #1~$#2{#4}$%
\fi}
\else%
{\ifmmode%
#2{#3}#1#2{#3+1}#1\dots #1#2{#4}%
\else%
$#2{#3}$#1~$#2{#3+1}$#1 \dots #1~$#2{#4}$%
\fi}
\fi}
\author{Martin Rubey}
\newcommand{\T}{\Set{T}}

\newcommand{\lmin}{\underline{l}}
\newcommand{\lmax}{\overline{l}}
\newcommand{\lind}{l}
\newcommand{\mmin}{\underline{m}}
\newcommand{\mmax}{\overline{m}}

\newcommand{\Lmin}{\underline{L}}
\newcommand{\Lmax}{\overline{L}}

\newcolumntype{S}{!{$\makebox[0pt]{$\mid$}$}}
\newcolumntype{D}{!{$\makebox[0pt]{$\lower2.4pt\vbox{\baselineskip=2.2pt\lineskiplimit=0pt\hbox{.}\hbox{.}\hbox{.}\hbox{.}\hbox{.}}$}$}}
\newcommand{\SText}[1]{\multicolumn{1}{Sc}{#1}}
\newcommand{\SP}{\multicolumn{1}{S!.}{}}
\newcommand{\DText}[1]{\multicolumn{1}{Dc}{#1}} 
\newcommand{\DP}{\multicolumn{1}{D!.}{}}

\begin{document}
\begin{abstract}
  We show that the $h$-vector of a ladder determinantal ring cogenerated by
  $M=[u_1\mid v_1]$ is log-concave. Thus we prove an instance of a conjecture
  of Stanley, resp. Conca and Herzog.
\end{abstract}
\title[]{The $h$-vector of a ladder determinantal ring cogenerated by $2\times
  2$ minors is log-concave} \maketitle
\thanks{In honour of Miriam Rubey, at the occasion of her second birthday}
\section{Introduction}
\begin{dfn}
  A sequence of real numbers $\seq{a_}{1}{n}$ is \Dfn{logarithmically concave},
  for short \Dfn{log-concave}, if $a_{i-1}a_{i+1}\le a_i^2$ for
  $i\in\{\seq{}{2}{n-1}\}$.
\end{dfn}
Numerous sequences arising in combinatorics and algebra have, or seem to have
this property. In the paper \cite{Stanley1989} written in 1989, Richard Stanley
collected various results on this topic. (For an update see \cite{Brenti1994}.)
There he also stated the following conjecture:
\begin{cnj}\label{cnj:Stanley}
  Let $R=R_0\oplus R_1\oplus\dots$ be a graded (Noetherian) Cohen-Macaulay (or
  perhaps Gorenstein) \emph{domain} over a field $K=R_0$, which is generated by
  $R_1$ and has Krull dimension $d$. Let $H(R,m)=\dim_K R_m$ be the Hilbert
  function of $R$ and write
  \begin{equation*}
    \sum_{m\ge0}H(R,m)x^m=(1-x)^{-d}\sum_{i=0}^s h_i x^i.
  \end{equation*}
  Then the sequence $\seq{h_}{0}{s}$ is log-concave.
\end{cnj}
The sequence $\seq{h_}{0}{s}$ is called the \Dfn{$h$-vector} of the ring.
Orginally the question was to decide whether a given sequence can arise as the
$h$-vector of some ring. In this sense the validity of the conjecture would
imply that log-concavity was a necessary condition on the $h$-vector.

It is now known however \cite{NiesiRobbiano1992,Brenti1994} that Stanley's
conjecture is not true in general. Several natural weakenings have been
considered, but are still open. For example, Aldo Conca and J\"urgen Herzog
conjectured that the $h$-vector would be log-concave for the special case where
$R$ is a ladder determinantal ring. (Note that ladder determinantal rings are
Cohen-Macaulay, as was shown in \cite[Corollary~4.10]{HerzogTrung1992}, but not
necessarily Gorenstein.) We will prove the conjecture of Conca and Herzog in the
simplest case, i.e., where $R$ is a ladder determinantal ring cogenerated by
$2\times 2$ minors, see Corollary~\ref{cor:main}.

In the case of ladder determinantal rings the $h$-vector has a nice
combinatorial interpretation. This follows from work of Abhyankar and Kulkarni
\cite{Abhyankar1988,AbhyankarKulkarni1989,Kulkarni1993,Kulkarni1996}, Bruns,
Conca, Herzog, and Trung
\cite{BrunsHerzog1992,Conca1995,ConcaHerzog1994,HerzogTrung1992}. In the
following paragraphs, which are taken almost verbatim from
\cite{KrattenthalerRubey2001}, we will explain these matters.
\section{Ladders, ladder determinantal rings and non-intersecting
  lattice paths}\label{s:Defs} 
First we have to introduce the notion of a ladder:
\begin{dfn}\label{dfn:ladder}
  Let $\Mat X=(x_{i,j})_{0\le i\le b, 0\le j\le a}$ be a $(b+1)\times (a+1)$
  matrix of indeterminates. Let $\Mat Y=(y_{i,j})_{0\le i\le b, 0\le j\le a}$
  be another matrix of the same dimensions, with the property that
  $y_{i,j}\in\{0,x_{i,j}\}$, and if $y_{i,j}=x_{i,j}$ and
  $y_{i^\prime,j^\prime}=x_{i^\prime,j^\prime}$, where $i\le i^\prime$ and
  $j\le j^\prime$ then $y_{r,s}=x_{r,s}$ for all $r$ and $s$ with $i\le r\le
  i^\prime$ and $j\le s\le j^\prime$. Such a matrix $\Mat Y$ is called a
  \Dfn{ladder}.
  
  A \Dfn{ladder region} $L$ is a subset of $\mathbb Z^2$ with the property that
  if $(i,j)$ and $(i^\prime,j^\prime)\in L$, $i\le i^\prime$ and $j\ge
  j^\prime$ then $(r,s)\in L$ for all $r\in\{\seq{}{i}{i^\prime}\}$ and
  $s\in\{\seq{}{j^\prime}{j}\}$.  Clearly, a ladder region can be described by
  two weakly increasing functions $\Lmin$ and $\Lmax$, such that $L$ is exactly
  the set of points $\{(i,j):\Lmin(i)\le j\le\Lmax(i)\}$.
  
  We associate with $\Mat Y$ a ladder region $L\subset\mathbb Z^2$ via
  $(j,b-i)\in L$ if and only if $y_{i,j}=x_{i,j}$.
\end{dfn}
In Figure~\ref{fig:latticepaths}.a an example of a ladder with $a=8$ and $b=9$
is shown, the corresponding ladder region is shown in
Figure~\ref{fig:latticepaths}.b.

Now we can define the ring we are dealing with:
\begin{dfn}\label{dfn:ring}
  Given a $(b+1)\times (a+1)$ matrix $\Mat Y$ which is a ladder, fix a
  ``bivector'' $M=[\seq{u_}{1}{n}\mid\seq{v_}{1}{n}]$ of integers with
  $1\le\seq[<]{u_}{1}{n}\le b+1$ and $1\le\seq[<]{v_}{1}{n}\le a+1$. By
  convention we set $u_{n+1}=b+2$ and $v_{n+1}=a+2$.
  
  Let $K[\Mat Y]$ denote the ring of all polynomials over some field $K$ in the
  $y_{i,j}$'s, where $0\le i\le b$ and $0\le j\le a$.  Furthermore, let
  $I_M(\Mat Y)$ be the ideal in $K[\Mat Y]$ that is generated by those $t\times
  t$ minors of $\Mat Y$ that contain only nonzero entries, whose rows form a
  subset of the last $u_t-1$ rows \emph{or} whose columns form a subset of the
  last $v_t-1$ columns, $t\in\{\seq{}{1}{n+1}\}$.  Thus, for $t=n+1$ the rows
  and columns of minors are unrestricted.
  
  The ideal $I_M(\Mat Y)$ is called a \Dfn{ladder determinantal ideal generated
    by the minors defined by $M$}. We call $R_M(\Mat Y)=K[\Mat Y]/I_M(\Mat Y)$
  the \Dfn{ladder determinantal ring cogenerated by the minors defined by $M$},
  or, in abuse of language, the \Dfn{ladder determinantal ring cogenerated by
    $M$}.
\end{dfn}
Note that we could restrict ourselves to the case $u_1=v_1=1$, because all the
elements of $\Mat Y$ that are in one of the last $u_1-1$ rows or in one of the
last $v_1-1$ columns are in the ideal. 

Next, we introduce the combinatorial objects that will accompany us
throughout the rest of this paper:
\begin{dfn}
  A \Dfn{two-rowed array of length $k$} is a pair of strictly increasing
  sequences of integers, both of length $k$. A two-rowed array $T=\left(
  \begin{smallmatrix}
    a_1&a_2&\dots&a_k\\
    b_1&b_2&\dots&b_k
  \end{smallmatrix}\right)$ is \Dfn{bounded} by $A=(A_1,A_2)$ and 
  $E=(E_1,E_2)$, if
  \begin{equation*}
    \begin{array}{lr!{\le}c!{\le}l}
                & A_1  &\seq[<]{a_}{1}{k}& E_1-1\\
      \text{and}& A_2+1&\seq[<]{b_}{1}{k}& E_2.
    \end{array}
  \end{equation*}
  Given \emph{any} subset $L$ of $\mathbb Z^2$, we say that the two-rowed array
  $T$ is \Dfn{in $L$}, if $(a_i,b_i)\in L$ for $i\in\{\seq{}{1}{k}\}$. By
  $\T^L_k(A\mapsto E)$ we will denote the set of two-rowed arrays of length
  $k$, bounded by $A$ and $E$ which are in $L$.  The \Dfn{total length} of a
  family of two-rowed arrays is just the sum of the lengths of its members.

  Let $T_1=\left(
    \begin{smallmatrix}
      a_1&a_2&\dots&a_k\\
      b_1&b_2&\dots&b_k
    \end{smallmatrix}\right)$ and $T_2=\left(
    \begin{smallmatrix}
      x_1&x_2&\dots&x_l\\
      y_1&y_2&\dots&y_l
    \end{smallmatrix}\right)$ be two-rowed arrays bounded by
  $A^{(1)}=(A^{(1)}_1,A^{(1)}_2)$ and $E^{(1)}=(E^{(1)}_1,E^{(1)}_2)$ and
  $A^{(2)}=(A^{(2)}_1,A^{(2)}_2)$ and $E^{(2)}=(E^{(2)}_1,E^{(2)}_2)$
  respectively. Set $a_{k+1}=E^{(1)}_1$ and $b_0=A^{(1)}_2$. We say that $T_1$
  and $T_2$ \Dfn{intersect} if there are indices $I$ and $J$ such that
  \begin{align*}\tag{$\times$}\label{e:intersect}
        x_J&\le a_I\\
    b_{I-1}&\le y_J
  \end{align*}
  where $1\le I\le k+1$ and $1\le J\le l$. A family of two-rowed arrays is
  \Dfn{non-intersecting} if no two arrays in it intersect.
\end{dfn}
Note that a two-rowed array in $\T^L_k(A\mapsto E)$ can be visualized by a
lattice path with east and north steps, that starts in $A$ and terminates in
$E$ and has exactly $k$ north-east turns which are all in $L$: Each pair
$(a_i,b_i)$ of a two-rowed array $\left(
  \begin{smallmatrix}
    a_1&a_2&\dots&a_k\\
    b_1&b_2&\dots&b_k
  \end{smallmatrix}\right)$ then corresponds to a north-east turn of the
lattice path. It is easy to see that Condition~\eqref{e:intersect} holds if and
only if the lattice paths corresponding to $T_1$ and $T_2$ intersect.

For an example see Figure~\ref{fig:latticepaths}.c, where the three two-rowed
arrays
\begin{equation*} T^{(1)}=
  \begin{pmatrix}
     2 & 3 \\
     6 & 7 
  \end{pmatrix}\text{, } T^{(2)}=
  \begin{pmatrix}
     3 & 5\\
     4 & 6
  \end{pmatrix}\text{ and } T^{(3)}=
  \begin{pmatrix}
    2 & 4 & 6\\
    1 & 3 & 4
  \end{pmatrix}
\end{equation*}
bounded by $A^{(1)}=(0,3)$, $A^{(2)}=(0,2)$, $A^{(3)}=(0,0)$ and $E^{(1)}=(5,9)$,
$E^{(2)}=(7,9)$, $E^{(3)}=(8,9)$ are shown as lattice paths. The points of the
ladder-region $L$ are drawn as small dots, the circles indicate the start- and
endpoints and the big dots indicate the north-east turns.
\def\x{\mbox{\tiny$\bullet$}}
\def\X{\mbox{\Large$\bullet$}} 
\def\O{\mbox{\Huge$\circ$}} 
\def\LadderMatrix{\POS
                                                (7,9)*{x_{0,7}},(8,9)*{x_{0,8}},
                                                (7,8)*{x_{1,7}},(8,8)*{x_{1,8}}, 
                                                (7,7)*{x_{2,7}},(8,7)*{x_{2,8}},
                (5,9)*{x_{0,5}},(6,9)*{x_{0,6}},
(4,8)*{x_{1,4}},(5,8)*{x_{1,5}},(6,8)*{x_{1,6}},
                                (2,7)*{x_{2,2}},(3,7)*{x_{2,3}},
(4,7)*{x_{2,4}},(5,7)*{x_{2,5}},(6,7)*{x_{2,6}},
                                (2,6)*{x_{3,2}},(3,6)*{x_{3,3}},
(4,6)*{x_{3,4}},(5,6)*{x_{3,5}},(6,6)*{x_{3,6}},
                                (2,5)*{x_{4,2}},(3,5)*{x_{4,3}},
(4,5)*{x_{4,4}},(5,5)*{x_{4,5}},(6,5)*{x_{4,6}},
                                (2,4)*{x_{5,2}},(3,4)*{x_{5,3}},
(4,4)*{x_{5,4}},(5,4)*{x_{5,5}},(6,4)*{x_{5,6}},
(0,3)*{x_{6,0}},(1,3)*{x_{6,1}},(2,3)*{x_{6,2}},(3,3)*{x_{6,3}},
(4,3)*{x_{6,4}},(5,3)*{x_{6,5}},(6,3)*{x_{6,6}},
(0,2)*{x_{7,0}},(1,2)*{x_{7,1}},(2,2)*{x_{7,2}},(3,2)*{x_{7,3}},
(4,2)*{x_{7,4}},(5,2)*{x_{7,5}},           
(0,1)*{x_{8,0}},(1,1)*{x_{8,1}},(2,1)*{x_{8,2}},(3,1)*{x_{8,3}},
(4,1)*{x_{8,4}},                          
(0,0)*{x_{9,0}},(1,0)*{x_{9,1}},(2,0)*{x_{9,2}},(3,0)*{x_{9,3}},
(4,0)*{x_{9,4}},
(0,4)*{0},(1,4)*{0},
(0,5)*{0},(1,5)*{0},
(0,6)*{0},(1,6)*{0},
(0,7)*{0},(1,7)*{0},
(0,8)*{0},(1,8)*{0},(2,8)*{0},(3,8)*{0},
(0,9)*{0},(1,9)*{0},(2,9)*{0},(3,9)*{0},(4,9)*{0},
(8,0)*{0},(7,0)*{0},(6,0)*{0},(5,0)*{0},
(8,1)*{0},(7,1)*{0},(6,1)*{0},(5,1)*{0},
(8,2)*{0},(7,2)*{0},(6,2)*{0},           
(8,3)*{0},(7,3)*{0},
(8,4)*{0},(7,4)*{0},
(8,5)*{0},(7,5)*{0},
(8,6)*{0},(7,6)*{0}}
\def\Ladder{\POS
                                                         (7,9)*{\x},(8,9)*{\x},
                                                         (7,8)*{\x},(8,8)*{\x}, 
                                                         (7,7)*{\x},(8,7)*{\x},
                                            (4,9)*{  },(5,9)*{\x},(6,9)*{\x}, 
                                            (4,8)*{\x},(5,8)*{\x},(6,8)*{\x},
                      (2,7)*{\x},(3,7)*{\x},(4,7)*{\x},(5,7)*{\x},(6,7)*{\x},
                      (2,6)*{\x},(3,6)*{\x},(4,6)*{\x},(5,6)*{\x},(6,6)*{\x},
                      (2,5)*{\x},(3,5)*{\x},(4,5)*{\x},(5,5)*{\x},(6,5)*{\x},
                      (2,4)*{\x},(3,4)*{\x},(4,4)*{\x},(5,4)*{\x},(6,4)*{\x},
(0,3)*{\x},(1,3)*{\x},(2,3)*{\x},(3,3)*{\x},(4,3)*{\x},(5,3)*{\x},(6,3)*{\x},
(0,2)*{\x},(1,2)*{\x},(2,2)*{\x},(3,2)*{\x},(4,2)*{\x},(5,2)*{\x},
(0,1)*{\x},(1,1)*{\x},(2,1)*{\x},(3,1)*{\x},(4,1)*{\x},
(0,0)*{\x},(1,0)*{\x},(2,0)*{\x},(3,0)*{\x},(4,0)*{\x},
(-1,0)*{0},(-1,1)*{1},(-1,2)*{2},(-1,3)*{3},(-1,4)*{4},(-1,5)*{5},(-1,6)*{6},
(-1,7)*{7},(-1,8)*{8},(-1,9)*{9},
(0,-1)*{0},(1,-1)*{1},(2,-1)*{2},(3,-1)*{3},(4,-1)*{4},(5,-1)*{5},(6,-1)*{6},
(7,-1)*{7},(8,-1)*{8}}
\def\StartEnd{\POS (0,3)*{\O},(0,2)*{\O},(0,0)*{\O},(5,9)*{\O},(7,9)*{\O},
                   (8,9)*{\O}}
\def\StartEndI{\POS (0,3)*i{\O},(0,2)*i{\O},(0,0)*i{\O},(5,9)*i{\O},(7,9)*i{\O},
                   (8,9)*i{\O}}
\def\NETurns{\POS (2,6)*{\X},(3,7)*{\X},(3,4)*{\X},(4,3)*{\X},(5,6)*{\X},
                  (6,4)*{\X},(2,1)*{\X}}
\def\PathNEA{\POS
@i @={(2,0),(0,3),(1,0),(0,1),(2,0),(0,2)},
(-0.08,3.08)="prev", @@{;p+"prev";"prev";**@{-}="prev"}}
%
\def\PathNEB{\POS
@i @={(3,0),(0,2),(2,0),(0,2),(2,0),(0,3)},
(-0.08,2.08)="prev", @@{;p+"prev";"prev";**@{-}="prev"}}
%
\def\PathNEC{\POS
@i @={(2,0),(0,1),(2,0),(0,2),(2,0),(0,1),(2,0),(0,5)},
(-0.08,0.08)="prev", @@{;p+"prev";"prev";**@{-}="prev"}}
\begin{figure}[h]
  \begin{tabular}{ccc}
\raisebox{-14pt}{\tiny\begin{xy}<12.5pt,0pt>:<0pt,11.45pt>::*+\xybox{\LadderMatrix}*\frm{-}\end{xy}}
&
\begin{xy}<10pt,0pt>:*+\xybox{\Ladder\StartEndI}*\frm{-}\end{xy}
&
\begin{xy}<10pt,0pt>:*+\xybox{\Ladder\StartEnd\NETurns
                               \PathNEA\PathNEB\PathNEC}*\frm{-}\end{xy}
\\\\
\parbox[t]{3.5cm}{a. a ladder with $a=8$ and $b=9$}
&                
\parbox[t]{3.5cm}{b. the corresponding ladder region}
&                
\parbox[t]{3.5cm}{c. a triple of non-intersecting lattice paths in this ladder}
  \end{tabular}
  \caption{\label{fig:latticepaths}}
\end{figure}
\section{A combinatorial interpretation of the $h$-vector of a ladder 
determinantal ring}
We are now ready to state the theorem which reveals the combinatorial nature of
the $h$-vector of $R_M(\Mat Y)=K[\Mat Y]/I_M(\Mat Y)$, the ladder determinantal
ring cogenerated by $M$.
\begin{thm}\label{thm:comb}
  Let $\Mat Y=(y_{i,j})_{0\le i\le b,\ 0\le j\le a}$ be a ladder and let
  $M=[\seq{u_}{1}{n}\mid\seq{v_}{1}{n}]$ be a bivector of integers with
  $1\le\seq[<]{u_}{1}{n}\le a+1$ and $1\le\seq[<]{v_}{1}{n}\le b+1$. For
  $i\in\{\seq{}{1}{n}\}$ let
  \begin{align*}
    A^{(i)}&=(0,u_{n+1-i}-1)\\
    E^{(i)}&=(a-v_{n+1-i}+1,b).  
    \intertext{Let $L^{(n)}=L$ be the ladder region associated
      with $\Mat Y$ and for $i\in\{\seq{}{1}{n-1}\}$ let}
    L^{(i)}&=\{(x,y)\in L^{(i+1)}:x\le E^{(i)}_1, y\ge A^{(i)}_2
    \text{ and } (x+1,y-1)\in L^{(i+1)}\}.
    \intertext{Finally, for $i\in\{\seq{}{1}{n}\}$ let}
    B^{(i)}&=\{(x,y)\in L^{(i)}:(x+1,y-1)\notin L^{(i)}\}.                                  
  \end{align*}
  and let  $d$ be the cardinality of $\bigcup_{i=1}^n B^{(i)}$.

  Then, under the assumption that all of the points $A^{(i)}$ and $E^{(i)}$,
  $i\in\{\seq{}{1}{n}\}$, lie inside the ladder region $L$, the Hilbert series
  of the ladder determinantal ring $R_M(\Mat Y)=K[\Mat Y]/I_M(\Mat Y)$ equals
  \begin{equation*}
    \sum _{\ell\ge0}\dim_K R_M(\Mat Y)_\ell\,z^\ell
    =\frac{\sum_{\ell\ge0}\size{\T^L_\ell(\Mat A\mapsto\Mat E)}z^\ell}
          {(1-z)^d}.
  \end{equation*}
  Here, $R_M(\Mat Y)_\ell$ denotes the homogeneous component of degree $\ell$
  in $R_M(\Mat Y)$ and $\size{\T^L_\ell(\Mat A\mapsto\Mat E)}$ is the number of
  non-intersecting families of two-rowed arrays with total length $\ell$, such
  that the $i$\textsuperscript{th} two-rowed array is bounded by $A^{(i)}$ and
  $E^{(i)}$ and is in $L^{(i)}\setminus B^{(i)}$ for $i\in\{\seq{}{1}{n}\}$.
\end{thm}
The sets $B^{(i)}$, $i\in\{\seq{}{1}{n}\}$ can be visualized as being the
lower-right boundary of $L^{(i)}$. Viewed as a path, there are exactly
$E^{(i)}_1-A^{(i)}_1+E^{(i)}_2-A^{(i)}_2+1$ lattice points on $B^{(i)}$, but
not all of them are necessarily in $L$. However, if $L$ is an upper ladder,
that is, $(a,0)\in L$, then this must be the case and we have
\begin{align*}
  d&=\sum_{i=1}^n\left(E^{(i)}_1-A^{(i)}_1+E^{(i)}_2-A^{(i)}_2+1\right)\\
   &=\sum_{i=1}^n\left(a-v_{n+1-i}+1 + b -u_{n+1-i}+1+1\right)\\
   &=n(a+b+3)-\sum_{i=1}^n\left(u_i+v_i\right),
\end{align*}
as in \cite{KrattenthalerRubey2001}.

In Figure~\ref{fig:LadderFace}.a, an example for a ladder region $L$ with $a=8$
and $b=9$ is given. The small dots represent elements of $L$, the circles on
the left and on the top of $L$ represent the points $A^{(i)}$ and $E^{(i)}$,
$i\in\{1,2,3\}$ that are specified by the minor $M=[1,3,4\mid 1,2,4]$. The
dotted lines indicate the lower boundary of $L^{(i)}$. Note that the point
$(4,9)$ is not an element of $L$. Therefore, in this example we have
\begin{equation*}
  d=n(a+b+3)-\sum_{i=1}^n\left(u_i+v_i\right)-1=44.
\end{equation*}
\def\Face{\POS (6,8)*{\X},(2,7)*{\X},(3,6)*{\X},(6,6)*{\X},(2,5)*{\X},
               (5,8)*{\X},(4,5)*{\X},(3,3)*{\X},(1,2)*{\X},(2,1)*{\X},(3,1)*{\X}}
\def\NETurns{\POS (2,7)*{\X},(3,6)*{\X},(4,5)*{\X},(2,1)*{\X}}
\def\LowestPaths{\POS
%
@i @={(4,0),(0,2),(1,0),(0,1),(1,0),(0,4),(2,0),(0,2)},
(0,0)="prev", @@{;p+"prev";"prev";**@{.}="prev"}
%
@i @={(3,0),(0,1),(1,0),(0,1),(1,0),(0,4),(2,0),(0,1)},
(0,2)="prev", @@{;p+"prev";"prev";**@{.}="prev"}
%
@i @={(2,0),(0,1),(1,0),(0,1),(1,0),(0,4),(1,0)},
(0,3)="prev", @@{;p+"prev";"prev";**@{.}="prev"}}
%
\def\PathA{\POS
@i @={(2,0),(0,4),(2,0),(0,2),(1,0)},
(-0.08,3.08)="prev", @@{;p+"prev";"prev";**@{-}="prev"}}
%
\def\PathNEA{\POS
@i @={(2,0),(0,4),(3,0),(0,2)},
(-0.08,3.08)="prev", @@{;p+"prev";"prev";**@{-}="prev"}}
%
\def\PathB{\POS 
@i @={(3,0),(0,4),(2,0),(0,2),(2,0),(0,1)},
(-0.08,2.08)="prev", @@{;p+"prev";"prev";**@{-}="prev"}}
%
\def\PathNEB{\POS
@i @={(3,0),(0,4),(4,0),(0,3)},
(-0.08,2.08)="prev", @@{;p+"prev";"prev";**@{-}="prev"}}
%
\def\PathC{\POS
@i @={(2,0),(0,1),(2,0),(0,4),(2,0),(0,2),(2,0),(0,2)},
(-0.08,0.08)="prev", @@{;p+"prev";"prev";**@{-}="prev"}}
%
\def\PathNEC{\POS
@i @={(2,0),(0,1),(2,0),(0,4),(4,0),(0,4)},
(-0.08,0.08)="prev", @@{;p+"prev";"prev";**@{-}="prev"}}
\begin{figure}[p]
  \begin{tabular}{cc}
\begin{xy}<11pt,0pt>:*+\xybox{\Ladder\StartEnd\LowestPaths}*\frm{-}\end{xy}
&
\begin{xy}<11pt,0pt>:*+\xybox{\Ladder\StartEnd\LowestPaths\Face}*\frm{-}\end{xy}
\\\\
\parbox[t]{4.5cm}{a. a ladder region with $a=8$ and $b=9$}
&
\parbox[t]{4.5cm}{b. a $10$ dimensional face of $\Delta_{[1,3,4\mid 1,2,4]}(\Mat
  Y)$}
  \end{tabular}
  \caption{\label{fig:LadderFace}}
\end{figure}
\begin{figure}[p]
  \begin{equation*}
\begin{xy}<11pt,0pt>:*+\xybox{\Ladder\StartEnd\LowestPaths\Face
                              \PathA}*\frm{-}\end{xy}
\quad                  
\begin{xy}<11pt,0pt>:*+\xybox{\Ladder\StartEnd\LowestPaths\Face
                              \PathA\PathB}*\frm{-}\end{xy}
\quad              
\begin{xy}<11pt,0pt>:*+\xybox{\Ladder\StartEnd\LowestPaths\Face
                              \PathA\PathB\PathC}*\frm{-}\end{xy}
  \end{equation*}
\caption{Constructing a family of non-intersecting lattice paths, such that the
  $i$\textsuperscript{th} path stays above $L^{(i)}$, $i\in\{1,2,3\}$}
  \label{fig:Paths}
\end{figure}
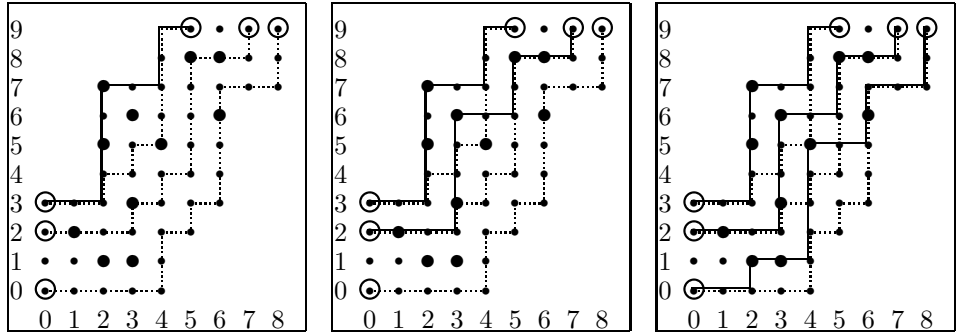
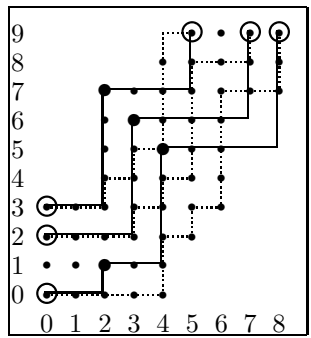
\begin{figure}[p]
  \begin{equation*}
\begin{xy}<11pt,0pt>:*+\xybox{\Ladder\StartEnd\LowestPaths\NETurns
                              \PathNEA\PathNEB\PathNEC}*\frm{-}\end{xy}
  \end{equation*}
  \caption{The corresponding family of non-intersecting lattice paths, 
           where the $i$\textsuperscript{th} path has north-east turns only in
           $L^{(i)}$ for $i\in\{1,2,3\}$}
  \label{fig:NEPaths}
\end{figure}
\begin{proof} 
  We will use results of J\"urgen Herzog and Ng\^o Vi\^et Trung. In Section~4
  of \cite{HerzogTrung1992}, ladder determinantal rings are introduced and
  investigated.
  
  We equip the indeterminates $x_{i,j}$, $i\in\{\seq{}{0}{b}\}$ and
  $j\in\{\seq{}{0}{a}\}$ with the following partial order:
  \begin{equation*}
      x_{i,j}\le x_{i^\prime,j^\prime}\text{ if }i\ge i^\prime
      \text{ and }j\le j^\prime.
  \end{equation*}
  A $t$-antichain in this partial order is a family of elements
  $x_{r_1,s_1},x_{r_2,s_2},\dots,x_{r_t,s_t}$ such that $\seq[<]{r_}{1}{t}$ and
  $\seq[<]{s_}{1}{t}$. Thus, a $t$-antichain corresponds to a sequence
  $(s_1,b-r_1),(s_2,b-r_2),\dots,(s_t,b-r_t)$ of $t$ points in the ladder
  region associated with $\Mat Y$, where each point lies strictly south-east of
  the previous ones.
  
  Let $D_t$ be the union of the last $u_t-1$ rows and the last $v_t-1$ columns
  of $\Mat Y$. Let $\Delta_M(\Mat Y)$ be the simplicial complex whose
  $k$-dimensional faces are subsets of elements of $\Mat Y$ of cardinality
  $k+1$ which do not contain a $t$-antichain in $D_t$ for
  $t\in\{\seq{}{1}{n+1}\}$. Let $f_k$ be the number of $k$-dimensional faces of
  $\Delta_M(\Mat Y)$ for $k\ge 0$. Then, Corollary~4.3 of
  \cite{HerzogTrung1992} states, that
  \begin{equation*}
    \dim_K R_M(\Mat Y)_\ell=\sum_{k\ge 0}\binom{\ell-1}{k}f_k.
  \end{equation*}

  In the following, we will find an expression for the numbers $f_k$ involving
  certain families of non-intersecting lattice paths.
  
  In Figure~\ref{fig:LadderFace}.b, a $10$-dimensional face of
  $\Delta_{[1,3,4\mid 1,2,4]}(\Mat Y)$ is shown, the elements of the face are
  indicated by bold dots. We will describe a modification of Viennot`s \lq
  light and shadow procedure\rq\ (with the sun in the top-left corner) that
  produces a family of $n$ non-intersecting lattice paths such that the
  $i$\textsuperscript{th} path runs from $A^{(i)}=(0,u_{n+1-i})$ to
  $E^{(i)}=(a-v_{n+1-i},b)$ and has north-east turns only in $L^{(i)}$, for
  $i\in\{\seq{}{1}{n}\}$.
  
  Imagine a sun in the top-left corner of the ladder region and a wall along
  the lower-right border $B^{(1)}$ of $L^{(1)}$. Then each lattice point
  $(r,s)$ that is either in $B^{(1)}$ or corresponds to an element $x_{s,b-r}$
  of the face casts a \lq shadow\rq\ $\{(x,y):x\ge r, y\le s\}$.
  
  The first path starts at $A^{(1)}$, goes along the north-east border of this
  shadow and terminates in $E^{(1)}$. In the left-most diagram of
  Figure~\ref{fig:Paths}, this is accomplished for the face shown in
  Figure~\ref{fig:LadderFace}.b.
  
  In the next step, we remove the wall on $B^{(1)}$ and all the elements of the
  face which correspond to lattice points lying on the first path. Then the
  procedure is iterated. See Figure~\ref{fig:Paths} for an example. Let $P$ be
  the resulting family of non-intersecting lattice paths.
  
  Now, for each $i\in\{\seq{}{1}{n}\}$, we remove all elements of the face
  except those which correspond to north-east turns of the
  $i$\textsuperscript{th} path and do not lie on $B^{(i)}$. In the example,
  $(5,8)$ is a north-east turn of the second path but lies on $B^{(2)}$,
  therefore the corresponding element $x_{1,5}$ of the face is removed. On the
  other hand, $(4,5)$ lies on $B^{(1)}$, but is a nort-east turn of the third
  path, so the corresponding element $x_{4,4}$ of the face is kept.
  
  This set of north-east turns defines another family of non-intersecting
  lattice paths $P^\prime$ that has the property that the
  $i$\textsuperscript{th} path has north-east turns only in $L^{(i)}$ for
  $i\in\{\seq{}{1}{n}\}$.
  
  We now want to count the number of faces of $\Delta_M(\Mat Y)$ that reduce
  under \lq light and shadow\rq\ to a given family of lattice paths $P^\prime$
  with this property. Clearly, $P^\prime$ can be translated into a family $P$
  of non-intersecting lattice paths such that the $i$\textsuperscript{th} path
  does not go below $B^{(i)}$ for $i\in\{\seq{}{1}{n}\}$. Note that the number
  of lattice points on such a family $P$ of paths is always equal to $d$,
  independently of the given face. Thus, if $m$ is the number of north-east
  turns of $P^\prime$, there are
  \begin{equation*}
    \binom{d-m}{k+1-m}
  \end{equation*}
  families of non-intersecting lattice paths $P$ that reduce to $P^\prime$.
  
  Hence, $f_k=\binom{d-m}{k+1-m}\size{\T^L_\ell(\Mat A\mapsto\Mat E)}$ and we
  obtain
  \begin{align*}
    \sum_{\ell\ge0}\dim_K R_M(\Mat Y)_\ell\,z^\ell
  &=\sum_{\ell\ge0}\bigg(\sum_{k\ge0}\binom{\ell-1}{k} f_k\bigg)z^\ell\\
  &=\sum_{\ell\ge0}\sum_{k\ge0}\binom{\ell-1}{k}
      \bigg(\sum_{m=0}^{k+1}\binom{d-m}{k+1-m}\size{\T^L_m(\Mat A\mapsto\Mat E)}\bigg)z^\ell\\
  &=\sum_{m\ge0}\size{\T^L_m(\Mat A\mapsto\Mat E)}\sum_{\ell\ge0}z^\ell
      \sum_{k\ge0}\binom{\ell-1}k\binom{d-m}{d-k-1},
  \end{align*}
  and if we sum the inner sum by means of the Vandermonde summation 
  (see for example \cite{GrahamKnuthPatashnik}, Section~5.1, (5.27)),
  \begin{align*}
    \sum_{\ell\ge0}\dim_K R_M(\Mat Y)_\ell\,z^\ell&=
    \sum_{m\ge0}\size{\T^L_m(\Mat A\mapsto\Mat E)}\sum_{\ell\ge0}z^\ell\binom{d+\ell-m-1}{d-1}\\
     &=\frac{\sum_{m\ge0}\size{\T^L_m(\Mat A\mapsto\Mat E)}z^m}{(1-z)^d}.
  \end{align*}
\end{proof}
\section{Log-concavity of the $h$-vector in the case $M=[u_1\mid v_1]$}
In this paper we will settle Stanley's conjecture when $R$ is a ladder
determinantal ring cogenerated by $M$, where $M$ is just a pair of integers,
i.e., $n=1$. We want to stress, however, that data strongly suggest that Conca
and Herzog's conjecture is also true for arbitrary $n$.

By the preceding theorem, in the case we are going to tackle, the sum
$\sum_{i=0}^s h_i x^i$ that appears in the conjecture is the
generating function $\sum_{k\ge0}\size{\T^L_k(A\mapsto E)}z^k$ of two-rowed
arrays bounded by $A$ and $E$ which are in the ladder region $L$.

As the bounds $A$ and $E$ will not be of any significance throughout the rest
of this paper, we will abbreviate $\T^L_k(A\mapsto E)$ to $\T^L_k$. We will
show that the $h$-vector is log-concave by constructing an injection from
$\T^L_{k+1}\times\T^L_{k-1}$ into $\T^L_k\times\T^L_k$. This injection will
involve some cut and paste operations that we now define:
\begin{dfn}
  Let $A$ and $X$ be two strictly increasing sequences of integers, such that
  the length of $X$ is the length of $A$ minus two, i.e.,
  $A=(a_1,a_2,\dots,a_{k+1})$ and $X=(x_1,x_2,\dots,x_{k-1})$ for some $k\ge
  1$. A \Dfn{cutting point of $A$ and $X$} is an index $l\in\{\seq{}{1}{k}\}$
  such that
  \begin{align*}\tag{$*$}\label{dfn:cut}
                        a_l&<x_l,\\
    \text{and}\quad x_{l-1}&<a_{l+1},
  \end{align*}  
  where we require the inequalities to be satisfied only if all variables are
  defined. Hence, $1$ is a cutting point if $a_1<x_1$, and $k$ is a cutting
  point if $x_{k-1}<a_{k+1}$.

  The \Dfn{image of $A$ and $X$ obtained by cutting at $l$} is
  \begin{equation*}
      \begin{array}{cccc@{\big\vert\!\!}c@{\!\!\big\vert}cccccc}
    a_1&a_2&\dots&\multicolumn{1}{c}{a_{l-1}}&a_l    &x_l    &x_{l+1} &\dots&x_{k-1}
\\[4pt]
    \hline\\[-8pt]
    x_1&x_2&\dots&x_{l-1}&\multicolumn{1}{c}{a_{l+1}}&a_{l+2}&\hdotsfor2&a_{k+1}
  \end{array}
  \end{equation*}
  Note that both the resulting sequences have length $k$.
\end{dfn}
\begin{lem}\label{l:possible}
  Let $A=(a_1,a_2,\dots,a_{k+1})$ and $X=(x_1,x_2,\dots,x_{k-1})$ be
  strictly increasing sequences of integers, such that the length of
  $X$ is the length of $A$ minus two. Then there exists at least one
  cutting point of $A$ and $X$.
\end{lem}
\begin{proof}
  If $a_l\ge x_l$ for $l\in\{\seq{}{1}{k-1}\}$ then
  $a_{k+1}>a_{k-1}\ge x_{k-1}$ and $k$ is a cutting point. Otherwise,
  let $l$ be minimal such that $a_l<x_l$. If $l=1$ then $1$ is a
  cutting point. Otherwise, because of the minimality of $l$, we have
  $a_{l+1}>a_{l-1}\ge x_{l-1}$, thus $l$ is a cutting point.
\end{proof}
\begin{dfn}
  Let $T=(T_1,T_2)\in\T_{k+1}\times\T_{k-1}$ be a pair of two-rowed arrays.
  Then a \Dfn{top cutting point of $T$} is a cutting point of the top rows of
  $T_1$ and $T_2$ and a \Dfn{bottom cutting point of $T$} is a cutting point of
  the bottom rows of $T_1$ and $T_2$.
  
  A pair $(l,m)$, where $l,m\in\{\seq{}{1}{k}\}$, such that $l$ is a top
  cutting point and $m$ is a bottom cutting point of $T_1$ and $T_2$ is a
  \Dfn{cutting point of $T$}. Cutting the top rows of $T$ at $l$ and the bottom
  rows at $m$ we obtain the \Dfn{image of $T$}.  Note that both of the
  two-rowed arrays in the image have length $k$. More pictorially, if $l<m$,
 \begin{equation*}
    \begin{array}{ccccccccccc}
      a_1&\hdotsfor2&a_l&x_l&\hdotsfor1&x_{m-1}&\hdotsfor2&x_{k-1}\\
      b_1&\hdotsfor3&b_{l+1}&\hdotsfor1&b_m&y_m&\dots&y_{k-1}\\[1pt]
      \hline
      x_1&\dots&x_{l-1}&a_{l+1}&\hdotsfor1&a_m&\hdotsfor3&a_{k+1}\\
      y_1&\hdotsfor2&y_l&\hdotsfor1&y_{m-1}&b_{m+1}&\hdotsfor2&b_{k+1},
    \end{array}
  \end{equation*}
  and similarly if $l\ge m$.
  
  For $T=(T_1,T_2)\in\T^L_{k+1}\times\T^L_{k-1}$, the pair $(l,m)$ is an
  \Dfn{allowed cutting point of $T$}, if both of the two-rowed arrays in the
  obtained image are in $L$.
\end{dfn}
In Lemma~\ref{l:allowed} we will prove that every pair of two-rowed arrays in
$\T^L_{k+1}\times\T^L_{k-1}$ has at least one allowed cutting point. This
motivates the following definition:
\begin{dfn}  
  Let $T=(T_1,T_2)\in\T^L_{k+1}\times\T^L_{k-1}$ a pair of two-rowed arrays as
  before. Consider all allowed cutting points $(\bar l,\bar m)$ of $T$. Select
  those with $\size{\bar l-\bar m}$ minimal. Among those, let $(l,m)$ be the
  pair which comes first in the lexicographic order. Then we call $(l,m)$ the
  \Dfn{optimal cutting point} of $T$.
\end{dfn}

Now we are ready to state our main theorem, which implies that Stanley's
conjecture is true, when $R$ is a ladder determinantal ring cogenerated by a
pair of integers $M$:
\begin{thm}\label{thm:main}
  Let $L$ be a ladder region. Let $T\in\T^L_{k+1}\times\T^L_{k-1}$. Define
  $I(T)$ to be the pair of two-rowed arrays obtained by cutting $T$ at its
  optimal cutting point. Then $I$ is well-defined and an injection from
  $\T^L_{k+1}\times\T^L_{k-1}$ into $\T^L_k\times\T^L_k$.
\end{thm}
\begin{cor}\label{cor:main}
  The $h$-vector of the ladder determinantal ring cogenerated by $M=[u_1\mid v_1]$ is
  log-concave.
\end{cor}
\begin{proof}[Proof of the corollary]
  By Theorem~\ref{thm:comb}, the $h$-vector of this ring is equal to the
  generating function $\sum_{k\ge0}\size{\T^L_k(A\mapsto E)}z^k$ of two-rowed
  arrays bounded by $A=(0,u_1-1)$ and $E=(a-v_1+1,b)$ which are in the ladder
  region $L$. By the preceding theorem, there is an injection from
  $\T^L_{k+1}(A\mapsto E)\times\T^L_{k-1}(A\mapsto E)$ into $\T^L_k(A\mapsto
  E)\times\T^L_k(A\mapsto E)$, thus
  \begin{equation*}
    \size{\T^L_{k+1}(A\mapsto E)}\cdot\size{\T^L_{k-1}(A\mapsto E)}
    \le\size{\T^L_k(A\mapsto E)}^2.
  \end{equation*}
\end{proof}

We will split the proof of Theorem~\ref{thm:main} in two parts. In
Section~\ref{s:welldef} we show that the mapping $I$ is well-defined, that is,
for any pair of two-rowed arrays $T\in\T^L_{k+1}\times\T^L_{k-1}$ there is an
allowed cutting point. Finally, in Section~\ref{s:inj}, we show that $I$ is
indeed an injection.

\section{The mapping $I$ is well-defined}\label{s:welldef}
\begin{lem}\label{l:allowed}
  Let $L$ be a ladder region. Then for every pair of two-rowed arrays in
  $\T^L_{k+1}\times\T^L_{k-1}$ there is an allowed cutting point $(l,m)$.
\end{lem}
For the proof of this lemma, we have to introduce some more notation: Let
$(T_1,T_2)\in\T^L_{k+1}\times\T^L_{k-1}$ with $T_1=\left(
\begin{smallmatrix}
  a_1&a_2&\dots&a_{k+1}\\
  b_1&b_2&\dots&b_{k+1}
\end{smallmatrix}\right)$ and $T_2=\left(
\begin{smallmatrix}
  x_1&x_2&\dots&x_{k-1}\\
  y_1&y_2&\dots&y_{k-1}
\end{smallmatrix}\right)$. 
We say that Inequality~\eqref{eq:topu} holds for an interval $[c,d]$
if
\begin{align}
  \label{eq:topu}\tag{$\overline{top}$}
  \Lmax(a_j)&\ge y_{j-1},\\
  \intertext{for $j\in[c,d]$. Inequality~\eqref{eq:topl} holds for an
    interval $[c,d]$ if}
  \label{eq:topl}\tag{\underline{$top$}}
  \Lmin(a_j)&\le y_{j-1},
\end{align}
for $j\in[c,d]$. Similarly, Inequality~\eqref{eq:bottomu} holds for an
interval $[c,d]$ if
\begin{align}
  \label{eq:bottomu}\tag{$\overline{bottom}$}
  \Lmax(x_{j-1})&\ge b_j,\\
  \intertext{for $j\in[c,d]$. Inequality~\eqref{eq:bottoml} holds for
  an interval $[c,d]$ if}
  \label{eq:bottoml}\tag{\underline{$bottom$}}
  \Lmin(x_{j-1})&\le b_j,
\end{align}
for $j\in[c,d]$, where $\Lmin$ and $\Lmax$ are as in
Definition~\ref{dfn:ladder}. We say that any of these inequalities holds for a
cutting point $(l,m)$ if it holds for the interval $[l+1,m]$ if $l<m$ and for
the interval $[m+1,l]$ if $m<l$.  Clearly, a cutting point $(l,m)$ is allowed
if and only if all of these inequalities hold for it.

Most of the work is done by the following lemma:
\begin{lem}\label{l:imps}
  Let $T=(T_1,T_2)\in\T^L_{k+1}\times\T^L_{k-1}$, $T_1=\left(
    \begin{smallmatrix}
      a_1&a_2&\dots&a_{k+1}\\
      b_1&b_2&\dots&b_{k+1}
    \end{smallmatrix}\right)$ and $T_2=\left(
    \begin{smallmatrix}
      x_1&x_2&\dots&x_{k-1}\\
      y_1&y_2&\dots&y_{k-1}
    \end{smallmatrix}\right)$. Let
  $\lmin$ and $\lmax$ be top cutting points, such that there is no top cutting
  point in the closed interval $[\lmin+1,\lmax-1]$. Similarly, let $\mmin$ and
  $\mmax$ be bottom cutting points, such that there is no bottom cutting point
  in the closed interval $[\mmin+1,\mmax-1]$. Then for both of the intervals
  $[\lmin+1,\lmax]$ and $[\mmin+1,\mmax]$,
  \begin{itemize}
  \item either \eqref{eq:topu} or \eqref{eq:bottomu} hold,
  \item either \eqref{eq:topl} or \eqref{eq:bottoml} hold,
  \item either \eqref{eq:topu} or \eqref{eq:topl} hold,
  \item either \eqref{eq:bottomu} or \eqref{eq:bottoml} hold.
  \end{itemize}
  Let $l_{min},l_{max},m_{min}$ and $m_{max}$ be the minimal and maximal top
  and bottom cutting points. Then we have
  \begin{itemize}
  \item \eqref{eq:topu} and \eqref{eq:bottoml} hold for
    $[2,\max(l_{min},m_{min})]$ and
  \item \eqref{eq:topl} and \eqref{eq:bottomu} hold for
    $[\min(l_{max},m_{max}),k]$.
  \end{itemize}
\end{lem}
\begin{proof}
  Suppose that \eqref{eq:topu} does not hold for the interval
  $[\lmin+1,\lmax]$. We claim that in that case there is an index
  $j\in[\lmin+1,\lmax-1]$ such that $a_j<x_j$: For, by hypothesis there is an
  index $i\in[\lmin+1,\lmax]$ such that $\Lmax(a_i)<y_{i-1}$. We have
  $\Lmax(a_i)<y_{i-1}\le\Lmax(x_{i-1})$ and because $\Lmax$ is a weakly
  increasing function, $a_i<x_{i-1}$. It follows that
  $a_{i-1}<a_i<x_{i-1}<x_i$. Thus, if $i=\lmax$ we choose $j=i-1$, otherwise
  $j=i$.
  
  The same statement is true if \eqref{eq:bottoml} does not hold for the
  interval $[\lmin+1,\lmax]$: In this case there must be an index
  $i\in[\lmin+1,\lmax]$ such that $\Lmin(x_{i-1})>b_i$. We conclude that
  $\Lmin(a_i)\le b_i<\Lmin(x_{i-1})$ and thus $a_i<x_{i-1}$.

  Next, we will use induction to prove that 
  \begin{align*}\tag{$**$}
    \label{eq:induction}
                    a_{\lind}  &<   x_{\lind}\\
    \text{and}\quad a_{\lind+1}&\le x_{\lind-1}
  \end{align*}
  for $\lind\in[\lmin+1,\lmax-1]$. We will first do an induction on $\lind$ to
  establish the claim for $\lind\in[j,\lmax-1]$.
  
  We start the induction at $\lind=j$: Above we already found that $a_j<x_j$.
  Therefore we must have $a_{j+1}\le x_{j-1}$, because otherwise $j$ would
  satisfy \eqref{dfn:cut} and hence were a top cutting point.
  
  Now suppose that \eqref{eq:induction} holds for a particular $\lind<\lmax-1$.
  Then $a_{\lind+1}\le x_{\lind-1}<x_{\lind+1}$, and, because there is no
  top cutting point at $\lind+1$, we have $a_{\lind+2}\le x_{\lind}$.
  
  Similarly, to establish \eqref{eq:induction} for $\lind\in[\lmin+1,j]$ we do
  a reverse induction on $\lind$. Suppose that \eqref{eq:induction} holds for a
  particular $\lind>\lmin+1$. Then $a_{\lind-1}<a_{\lind+1}\le x_{\lind-1}$,
  and, because there is no top cutting point at $\lind-1$, we have
  $a_{\lind}\le x_{\lind-2}$.
  
  Thus we obtain
  \begin{align*}
    &\Lmax(x_{\lmax-1})\ge\Lmax(x_{\lmax-2})\ge\Lmax(a_{\lmax})\ge b_{\lmax}
    \text{, and}\\
    &\Lmax(x_{\lind-1})\ge\Lmax(a_{\lind+1})\ge b_{\lind+1}    \ge b_{\lind},  
  \end{align*}
  which means that \eqref{eq:bottomu} holds for the interval $[\lmin+1,\lmax]$.

  Furthermore, 
  \begin{align*}
    &\Lmin(a_{\lmin+1})\le\Lmin(a_{\lmin+2})\le\Lmin(x_{\lmin})\le y_{\lmin}
    \text{, and}\\
    &\Lmin(a_{\lind})  \le\Lmin(x_{\lind-2})\le y_{\lind-2}    \le y_{\lind-1},  
  \end{align*}
  which means that \eqref{eq:topl} holds for the interval $[\lmin+1,\lmax]$.
  
  Next we show that \eqref{eq:topu} and \eqref{eq:bottoml} hold for the
  interval $[2,l_{min}]$: Assume that either of these inequalities does not
  hold for the interval $[2,l_{min}]$ and that $[2,l_{min}]$ does not contain a
  top cutting point except $l_{min}$. Then the above reverse induction implies
  that $a_1\le a_3<x_1$, which means that $1$ is a top cutting point. Thus,
  $l_{min}=1$ and the interval $[2,l_{min}]$ is empty.

  The other assertions are shown in a completely analogous fashion.
\end{proof}

We are now ready to establish Lemma~\ref{l:allowed}:
\begin{proof}[Proof of Lemma~\ref{l:allowed}]
  Let $T=(T_1,T_2)\in\T^L_{k+1}\times\T^L_{k-1}$. By Lemma~\ref{l:possible}
  there is at least one cutting point $(l,m)$ of $T$.  Let
  $l_{min},l_{max},m_{min}$ and $m_{max}$ be the minimal and maximal top and
  bottom cutting points of $T$ as before.
  
  If there is an index $j$ which is a top \emph{and} a bottom cutting point of
  $T$, then -- trivially -- $(j,j)$ is an allowed cutting point. Otherwise, we
  have to show that there is a cutting point $(l,m)$ for which \eqref{eq:topu},
  \eqref{eq:topl}, \eqref{eq:bottomu}, and \eqref{eq:bottoml} hold. Suppose
  that this is not the case.
  
  For the inductive proof which follows, we have to introduce a
  convenient indexing scheme for the sequence of top and bottom
  cutting points. Let
  \begin{alignat*}{2}
    m_{1,0}&=\max\{m:m<l_{min}\text{ and $m$ is a bottom
      cutting point}\},\\
    m_{i,0}&=\max\{m:m<l_{i-1,1}\text{ and $m$ is a bottom
      cutting point}\}&\quad&\text{for $i>1$,}\\
    \text{and}\quad l_{i,0}&=\max\{l:l<m_{i,1}\text{ and $l$ is a
      top cutting point}\} &\quad&\text{for $i\ge 1$,}
  \end{alignat*}
  where $m_{i,j+1}$ is the bottom cutting point directly after
  $m_{i,j}$, and $l_{i,j+1}$ is the top cutting point directly after
  $l_{i,j}$. Furthermore, we set $l_{0,1}=l_{min}$.
  
  More pictorially, we have the following sequence of top and bottom
  cutting points for $i\ge 1$:
  \begin{equation*}
    \dots<m_{i,0}<l_{i-1,1}<l_{i-1,2}<\dots<l_{i,0}<m_{i,1}<m_{i,2}
   <\dots<m_{i+1,0}<\cdots
  \end{equation*}
  If $m_{min}>l_{min}$, then $m_{1,0}$ does not exist, of course. Note
  that there are no bottom cutting points between $l_{i,1}$ and
  $l_{i+1,0}$, and there are no top cutting points between $m_{i,1}$
  and $m_{i+1,0}$.
  
  Suppose first that $m_{min}<l_{min}$. By induction on $i$, we will
  show that \eqref{eq:topu} and \eqref{eq:bottoml} hold for the
  cutting points $(l_{i-1,1},m_{i,0})$, where $i\ge 1$.  By
  Lemma~\ref{l:imps} we know that \eqref{eq:topu} and
  \eqref{eq:bottoml} are satisfied for the cutting point
  $(l_{min},m_{1,0})$, because $[m_{1,0}+1,l_{min}]\subseteq
  [2,l_{min}]$. It remains to perform the induction step, which we
  will divide into five simple steps.\medskip
  
  \emph{Step 1. \eqref{eq:topu} and \eqref{eq:bottoml} hold for the
    interval $[m_{i,0}+1,l_{i-1,1}]$.} This is just a restatement of
  the induction hypothesis, i.e., that \eqref{eq:topu} and
  \eqref{eq:bottoml} hold for the cutting point
  $(l_{i-1,1},m_{i,0})$.\medskip
  
  \emph{Step 2. Either \eqref{eq:bottomu} or \eqref{eq:topl} does not
    hold for the interval $[m_{i,0}+1,m_{i,1}]$.} Because of Step~1,
  not both of \eqref{eq:bottomu} and \eqref{eq:topl} can hold for
  $(l_{i-1,1},m_{i,0})$, lest this was an allowed cutting point. Thus
  either \eqref{eq:bottomu} or \eqref{eq:topl} does not hold for
  $[m_{i,0}+1,l_{i-1,0}+1]$. This interval is contained in
  $[m_{i,0}+1,m_{i,1}]$, thus the inequalities \eqref{eq:bottomu} and
  \eqref{eq:topl} cannot hold on this interval either.\medskip
  
  \emph{Step 3. \eqref{eq:topu} and \eqref{eq:bottoml} hold for
    $[l_{i,0}+1,m_{i,1}]$.} Suppose that \eqref{eq:bottomu} does not hold for
  $[m_{i,0}+1,m_{i,1}]$. Then, by Lemma~\ref{l:imps} we obtain that
  \eqref{eq:topu} and \eqref{eq:bottoml} hold for $[m_{i,0}+1,m_{i,1}]$,
  because this interval contains no bottom cutting points except $m_{i,1}$.
  The same is true, if \eqref{eq:topl} does not hold for
  $[m_{i,0}+1,m_{i,1}]$. Because $[l_{i,0}+1,m_{i,1}]$ is a subset of this
  interval, \eqref{eq:topu} and \eqref{eq:bottoml} hold for the cutting point
  $(l_{i,0},m_{i,1})$, or, equivalently, for the interval
  $[l_{i,0}+1,m_{i,1}]$.\medskip
  
  \emph{Step 4. Either \eqref{eq:bottomu} or \eqref{eq:topl} does not hold for
    $[l_{i,0}+1,l_{i,1}]$.} Because of Step~3, not both of \eqref{eq:bottomu}
  and \eqref{eq:topl} can hold for the cutting point $(l_{i,0},m_{i,1})$, nor
  for the greater interval $[l_{i,0}+1,l_{i,1}]$.\medskip
  
  \emph{Step 5. \eqref{eq:topu} and \eqref{eq:bottoml} hold for
    $[m_{i+1,0}+1,l_{i,1}]$.} The interval $[l_{i,0}+1,l_{i,1}]$ does not
  contain a top cutting point except $l_{i,1}$, thus by Lemma~\ref{l:imps}
  and Step~4 we see that \eqref{eq:topu} and \eqref{eq:bottoml} hold. Finally,
  because $[m_{i+1,0}+1,l_{i,1}]\subset [l_{i,0}+1,l_{i,1}]$,
  \eqref{eq:topu} and \eqref{eq:bottoml} hold for the cutting point
  $(l_{i,1},m_{i+1,0})$.\medskip
    
  If $l_{max}>m_{max}$, then we encounter a contradiction: Let $r$ be such that
  $m_{r,0}=m_{max}$. We have just shown that \eqref{eq:topu} and
  \eqref{eq:bottoml} hold for the cutting point $(l_{r-1,1},m_{r,0})$.
  Furthermore, by Lemma~\ref{l:imps}, \eqref{eq:bottomu} and \eqref{eq:topl}
  hold for $[m_{r,0},k]$ and thus also for $(l_{r-1,1},m_{r,0})$. Hence, this
  would be an allowed cutting point, contradicting our hypothesis.
  
  If $l_{max}<m_{max}$, let $r$ be such that $l_{r,0}=l_{max}$. By the
  induction (Step~3) we find that \eqref{eq:topu} and
  \eqref{eq:bottoml} hold for the cutting point $(l_{r,0},m_{r,1})$.
  Again, because of Lemma~\ref{l:imps}, we know that
  \eqref{eq:bottomu} and \eqref{eq:topl} holds for $[l_{r,0},k]$ and
  thus also for $(l_{r,0},m_{r,1})$. Hence, we had an allowed cutting
  point in this case also.
  
  The case that $m_1>l_1$ is completely analogous.
\end{proof}

\section{The mapping $I$ is an injection}\label{s:inj}
\begin{lem}
  The mapping $I$ defined above is an injection.
\end{lem}
\begin{proof} 
  Suppose that $I(T)=I(T^\prime)$ for $T=(T_1,T_2)$ and
  $T^\prime=(T_1^\prime,T_2^\prime)$, such that $T$ and $T^\prime$ are elements
  of $\T^L_{k+1}\times\T^L_{k-1}$. Let $(l,m)$ be the optimal cutting point of
  $T$, and let $(l^\prime,m^\prime)$ be the optimal cutting point of
  $T^\prime$.
  
  Observe that we can assume $\min(l,m,l^\prime,m^\prime)=1$, because
  the elements of $T$ and $T^\prime$ with index less than or equal to
  this minimum retain their position in $I(T)$. Likewise, we can
  assume that $\max(l,m,l^\prime,m^\prime)=k$.
  
  Furthermore, we can assume that $l\le l^\prime$, otherwise we
  exchange the meaning of $T$ and $T^\prime$. Thus, we have to
  consider the following twelve situations:
  \begin{equation*}
    \begin{array}{>(r<{)\quad} @{1=\:} c !{\le} c !{\le} c !{\le} c @{\:=k}}
      1 & l & l^\prime & m        & m^\prime\\
      2 & l & l^\prime & m^\prime & m\\  
      3 & l & m        & l^\prime & m^\prime\\
      4 & l & m        & m^\prime & l^\prime\\
      5 & l & m^\prime & l^\prime & m\\
      6 & l & m^\prime & m        & l^\prime\\
      7 & m & l        & l^\prime & m^\prime\\  
      8 & m & l        & m^\prime & l^\prime\\  
      9 & m & m^\prime & l        & l^\prime\\  
     10 & m^\prime & l        & l^\prime & m\\  
     11 & m^\prime & l        & m        & l^\prime\\  
     12 & m^\prime & m        & l        & l^\prime
    \end{array}
  \end{equation*}
  We shall divide these twelve cases into two portions according to whether
  $l\le m$ or not.

  \vbox{
  \subsection*{$\mathbf{\text{A: }l\le m}$}
  In the Cases~$(1)$--$(6)$, $(10)$ and $(11)$ we have $l\le m$, thus
  the pair of two-rowed arrays
  $T=(T_1,T_2)\in\T^L_{k+1}\times\T^L_{k-1}$ looks like
  \begin{equation*}
    \begin{array}{cccccccccc}
      a_1&\hdotsfor2&a_l&\SText{a_{l+1}}&\hdotsfor4&a_{k+1}\\
      b_1&\hdotsfor5&b_m&\SText{b_{m+1}}&\hdotsfor1&b_{k+1}\\[1pt]
      \hline
      x_1&\dots&x_{l-1}&\SText{x_l}&\hdotsfor4&x_{k-1}\\
      y_1&\hdotsfor4&y_{m-1}&\SText{y_m}&\hdotsfor1&y_{k-1}.
    \end{array}
  \end{equation*}
  Cutting at $(l,m)$ we obtain $I(T)\in\T^L_k\times\T^L_k$:
  \begin{equation*}
    \begin{array}{cccccccccc}
      a_1&\hdotsfor2&a_l&\SText{x_l}&\hdotsfor4&x_{k-1}\\
      b_1&\hdotsfor5&b_m&\SText{y_m}&\dots&y_{k-1}\\[1pt]
      \hline
      x_1&\dots&x_{l-1}&\SText{a_{l+1}}&\hdotsfor5&a_{k+1}\\
      y_1&\hdotsfor4&y_{m-1}&\SText{b_{m+1}}&\hdotsfor2&b_{k+1}.
    \end{array}
  \end{equation*}
  If $l=1$, then the top row of the second array in $I(T)$ is
  $(\seq{a_}{2}{k+1})$, if $m=k$, then the bottom row of the first array
  in $I(T)$ is $(\seq{b_}{1}{k})$.}

  \vbox{
  \subsection*{$\mathbf{\text{Case }(1)\text{, }
      1=l\le l^\prime\le m\le m^\prime=k}$} Given that $I(T)=I(T^\prime)$, the
  pair $T^\prime$ can be expressed in terms of the entries of $T$ as follows:
  \begin{equation*}\tag{$T^\prime$}
    \begin{array}{ccccccccccc}
    a_1
    &\SText{x_1}&\hdotsfor1&x_{l^\prime-1}&\DText{a_{l^\prime+1}}&\hdotsfor5
        &a_{k+1}\\
    b_1&\hdotsfor5&b_m
    &\SText{y_m}&\hdotsfor1&y_{k-1}&\DText{b_{k+1}}\\[1pt]
    \hline
    \SText{a_2}&\hdotsfor1&a_{l^\prime}&\DText{x_{l^\prime}}&\hdotsfor4
        &x_{k-1}&\\
    y_1&\hdotsfor4&y_{m-1}
    &\SText{b_{m+1}}&\hdotsfor1&b_k&\DP
  \end{array}
  \end{equation*}
  The vertical dots indicate the cut $(l^\prime,m^\prime)$ which would result
  in $I(T^\prime)$. We show that the cutting point $(l,m)=(1,m)$, indicated
  above by the vertical lines, is in fact an allowed cutting point for
  $T^\prime$: Cutting at $(1,m)$ yields
  \begin{equation}\tag{$\tilde T$}\label{eq:case1}
    \begin{array}{ccccccccccc}
    a_1
    &\SText{a_2}&\hdotsfor1&a_{l^\prime}&\DText{x_{l^\prime}}&\hdotsfor1
        &x_{m-1}&\hdotsfor2&x_{k-1}&\\
    b_1&\hdotsfor5&b_m
    &\SText{b_{m+1}}&\hdotsfor1&b_k&\DText{}\\[1pt]
    \hline
    \SText{x_1}&\hdotsfor1&x_{l^\prime-1}&\DText{a_{l^\prime+1}}&\hdotsfor1
        &a_m&\hdotsfor3&a_{k+1}\\
    y_1&\hdotsfor4&y_{m-1}
    &\SText{y_m}&\hdotsfor1&y_{k-1}&\DText{b_{k+1}.}
  \end{array}
  \end{equation}
  Note, that this is the same pair of two-rowed arrays we obtain by cutting $T$
  at $(l^\prime,m^\prime)$. We have to check that the pair of two-rowed arrays
  \eqref{eq:case1} is in the ladder region. 
  
  Clearly,
  \begin{equation*}
  (a_2,b_2),(a_3,b_3),\dots,(a_{l^\prime},b_{l^\prime})\text{ and }
  (x_1,y_1),(a_2,b_2),\dots,(x_{l^\prime-1},y_{l^\prime-1})
  \end{equation*}
    are in the ladder region, because these pairs appear also in $T$.
  Furthermore, the pairs
  \begin{align*}
    &(x_{l^\prime},b_{l^\prime+1}),(x_{l^\prime+1},b_{l^\prime+2}),\dots
     (x_{m-1},b_m)\\
   \text{and }
    &(a_{l^\prime+1},y_{l^\prime}),(a_{l^\prime+2},y_{l^\prime+1}),\dots
     (a_m,y_{m-1})
  \end{align*}
  appear in $I(T)$ and are therefore in the ladder region, too. All the other
  pairs, i.e.,
  \begin{gather*}
    (a_1,b_1)\text{ and }(x_m,b_{m+1}),(x_{m+1},b_{m+2}),\dots,(x_{k-1},b_k),\\
    (a_{m+1},y_m),(a_{m+2},y_{m+1}),\dots,(a_k,y_{k-1})\text{ and }(a_{k+1},b_{k+1}),
  \end{gather*}
  are unaffected by the cut and appear in $T^\prime$.  

  Thus we have that $(l,m)$ and $(l^\prime,m^\prime)$ are allowed cuts for $T$
  and $T^\prime$. We required that $(l,m)$ is optimal for $T$ and that
  $(l^\prime,m^\prime)$ is optimal for $T^\prime$, therefore we must have
  $l=l^\prime$ and $m=m^\prime$.}

In all the other cases the reasoning is very similar. Thus we only print the
pairs of two-rowed arrays $T^\prime$ and $\tilde T$ and leave it to the
reader to check that $\tilde T$ is in the ladder region.

  \vbox{
  \subsection*{$\mathbf{\text{Case }(2)\text{, }
      1=l\le l^\prime\le m^\prime\le m=k}$} The pair $T^\prime$ can be
  expressed in terms of the entries of $T$ as follows:
  \begin{equation*}\tag{$T^\prime$}
    \begin{array}{ccccccccccc}
    a_1
    &\SText{x_1}&\hdotsfor1&x_{l^\prime-1}&\DText{a_{l^\prime+1}}
        &\hdotsfor5&a_{k+1}\\
    b_1&\hdotsfor5&b_{m^\prime}&\DText{y_{m^\prime}}&\hdotsfor1&y_{k-1}
    &\SText{b_{k+1}}\\[1pt]
    \hline
    \SText{a_2}&\hdotsfor1&a_{l^\prime}&\DText{x_{l^\prime}}
        &\hdotsfor4&x_{k-1}&\\
    y_1&\hdotsfor4&y_{m^\prime-1}&\DText{b_{m^\prime+1}}&\hdotsfor1&b_k
    &\SP
  \end{array}
  \end{equation*}
  Cutting at $(l,m)$ yields  
  \begin{equation*}\tag{$\tilde T$}
    \begin{array}{ccccccccccc}
    a_1
    &\SText{a_2}&\hdotsfor1&a_{l^\prime}&\DText{x_{l^\prime}}&\hdotsfor1
        &x_{m^\prime-1}&\hdotsfor2&x_{k-1}&\\
    b_1&\hdotsfor5&b_{m^\prime}&\DText{y_{m^\prime}}&\hdotsfor1&y_{k-1}
    &\SText{}\\[1pt]
    \hline
    \SText{x_1}&\hdotsfor1&x_{l^\prime-1}&\DText{a_{l^\prime+1}}&\hdotsfor1
        &a_{m^\prime}&\hdotsfor3&a_{k+1}\\
    y_1&\hdotsfor4&y_{m^\prime-1}&\DText{b_{m^\prime+1}}&\hdotsfor1&b_k
    &\SText{b_{k+1}.}
    \end{array}
  \end{equation*}}

  \vbox{  
  \subsection*{$\mathbf{\text{Case }(3)\text{, }
      1=l\le m\le l^\prime\le m^\prime=k}$} The pair $T^\prime$ can be
  expressed in terms of the entries of $T$ as follows:
  \begin{equation*}\tag{$T^\prime$}
    \begin{array}{ccccccccccc}
    a_1
    &\SText{x_1}&\hdotsfor4&x_{l^\prime-1}&\DText{a_{l^\prime+1}}&\hdotsfor2
        &a_{k+1}\\
    b_1&\hdotsfor2&b_m
    &\SText{y_m}&\hdotsfor4&y_{k-1}&\DText{b_{k+1}}\\[1pt]
    \hline
    \SText{a_2}&\hdotsfor4&a_{l^\prime}&\DText{x_{l^\prime}}&\hdotsfor1&x_{k-1}
        &\\
    y_1&\hdotsfor1 &y_{m-1}
    &\SText{b_{m+1}}&\hdotsfor4&b_k&\DP
    \end{array}
  \end{equation*}
  Cutting at $(l,m)$ yields  
  \begin{equation*}\tag{$\tilde T$}
    \begin{array}{ccccccccccc}
    a_1
    &\SText{a_2}&\hdotsfor4&a_{l^\prime}&\DText{x_{l^\prime}}&\hdotsfor1
        &x_{k-1}&\\
    b_1&\hdotsfor2&b_m
    &\SText{b_{m+1}}&\hdotsfor4&b_k\\[1pt]
    \hline
    \SText{x_1}&\hdotsfor4&x_{l^\prime-1}&\DText{a_{l^\prime+1}}&\hdotsfor2
        &a_{k+1}\\
    y_1&\hdotsfor1 &y_{m-1}
    &\SText{y_m}&\hdotsfor4&y_{k-1}&\DText{b_{k+1}.}
  \end{array}
  \end{equation*}}

  \vbox{
  \subsection*{$\mathbf{\text{Case }(4)\text{, }
      1=l\le m\le m^\prime\le l^\prime=k}$} The pair $T^\prime$ can be
  expressed in terms of the entries of $T$ as follows:
  \begin{equation*}\tag{$T^\prime$}
    \begin{array}{ccccccccccc}
    a_1
    &\SText{x_1}&\hdotsfor7&x_{k-1}&\DText{a_{k+1}}\\
    b_1&\hdotsfor2&b_m
    &\SText{y_m}&\hdotsfor1&y_{m^\prime-1}&\DText{b_{m^\prime+1}}&\hdotsfor2
        &b_{k+1}\\[1pt]
    \hline
    \SText{a_2}&\hdotsfor7&a_k&\DText{}\\
    y_1&\hdotsfor1&y_{m-1}
    &\SText{b_{m+1}}&\hdotsfor1&b_{m^\prime}&\DText{y_{m^\prime}}&\hdotsfor1
        &y_{k-1}.
    \end{array}
  \end{equation*}
  Cutting at $(l,m)$ yields  
  \begin{equation*}\tag{$\tilde T$}
    \begin{array}{ccccccccccc}
    a_1&\SText{a_2}&\hdotsfor7&a_k&\DText{}\\
    b_1&\hdotsfor2&b_m
    &\SText{b_{m+1}}&\hdotsfor1&b_{m^\prime}&\DText{y_{m^\prime}}&\hdotsfor1
        &y_{k-1}\\[1pt]
    \hline
    \SText{x_1}&\hdotsfor7&x_{k-1}&\DText{a_{k+1}}\\
    y_1&\hdotsfor1&y_{m-1}
    &\SText{y_m}&\hdotsfor1&y_{m^\prime-1}&\DText{b_{m^\prime+1}}&\hdotsfor2
        &b_{k+1}.
    \end{array}
  \end{equation*}}

  \vbox{
  \subsection*{$\mathbf{\text{Case }(5)\text{, }
      1=l\le m^\prime\le l^\prime\le m=k}$} The pair $T^\prime$ can be
  expressed in terms of the entries of $T$ as follows:
  \begin{equation*}\tag{$T^\prime$}
    \begin{array}{ccccccccccc}
    a_1
    &\SText{x_1}&\hdotsfor4&x_{l^\prime-1}&\DText{a_{l^\prime+1}}&\hdotsfor2
        &a_{k+1}\\
    b_1&\hdotsfor2&b_{m^\prime}&\DText{y_{m^\prime}}&\hdotsfor4&y_{k-1}
    &\SText{b_{k+1}}\\[1pt]
    \hline
    \SText{a_2}&\hdotsfor4&a_{l^\prime}&\DText{x_{l^\prime}}&\hdotsfor1
        &x_{k-1}&\\
    y_1&\hdotsfor1&y_{m^\prime-1}&\DText{b_{m^\prime+1}}&\hdotsfor4&b_k
    &\SP
  \end{array}
  \end{equation*}
  Cutting at $(l,m)$ yields  
  \begin{equation*}\tag{$\tilde T$}
    \begin{array}{ccccccccccc}
    a_1
    &\SText{a_2}&\hdotsfor4&a_{l^\prime}&\DText{x_{l^\prime}}&\hdotsfor1
        &x_{k-1}&\\
    b_1&\hdotsfor2&b_{m^\prime}&\DText{y_{m^\prime}}&\hdotsfor4&y_{k-1}
    &\SText{}\\[1pt]
    \hline
    \SText{x_1}&\hdotsfor4&x_{l^\prime-1}&\DText{a_{l^\prime+1}}&\hdotsfor2
        &a_{k+1}\\
    y_1&\hdotsfor1&y_{m^\prime-1}&\DText{b_{m^\prime+1}}&\hdotsfor4&b_k
    &\SText{b_{k+1}.}
    \end{array}
  \end{equation*}}

  \vbox{
  \subsection*{$\mathbf{\text{Case }(6)\text{, }
      1=l\le m^\prime\le m\le l^\prime=k}$} The pair $T^\prime$ can be
  expressed in terms of the entries of $T$ as follows:
  \begin{equation*}\tag{$T^\prime$}
    \begin{array}{ccccccccccc}
    a_1&\SText{x_1}&\hdotsfor7&x_{k-1}&\DText{a_{k+1}}\\
    b_1&\hdotsfor2&b_{m^\prime}&\DText{y_{m^\prime}}&\hdotsfor1&y_{m-1}    
    &\SText{b_{m+1}}&\hdotsfor2&b_{k+1}\\[1pt]
    \hline
    \SText{a_2}&\hdotsfor7&a_k&\DText{}\\
    y_1&\hdotsfor1&y_{m^\prime-1}&\DText{b_{m^\prime+1}}&\hdotsfor1&b_m
    &\SText{y_m}&\hdotsfor1&y_{k-1}.
    \end{array}
  \end{equation*}
  Cutting at $(l,m)$ yields  
  \begin{equation*}\tag{$\tilde T$}
    \begin{array}{ccccccccccc}
    a_1&\SText{a_2}&\hdotsfor7&a_k&\DText{}\\
    b_1&\hdotsfor2&b_{m^\prime}&\DText{y_{m^\prime}}&\hdotsfor1&y_{m-1}
    &\SText{y_m}&\hdotsfor1&y_{k-1}\\[1pt]
    \hline
    \SText{x_1}&\hdotsfor7&x_{k-1}&\DText{a_{k+1}}\\
    y_1&\hdotsfor1&y_{m^\prime-1}&\DText{b_{m^\prime+1}}&\hdotsfor1&b_m    
    &\SText{b_{m+1}}&\hdotsfor2&b_{k+1}.
    \end{array}
  \end{equation*}}

  \vbox{
  \subsection*{$\mathbf{\text{Case }(10)\text{, }
      1=m^\prime\le l\le l^\prime\le m=k}$} The pair $T^\prime$ can be
  expressed in terms of the entries of $T$ as follows:
  \begin{equation*}\tag{$T^\prime$}
    \begin{array}{ccccccccccc}
    a_1&\hdotsfor2&a_l
    &\SText{x_l}&\hdotsfor1&x_{l^\prime-1}&\DText{a_{l^\prime+1}}&\hdotsfor2
        &a_{k+1}\\
    b_1
    &\DText{y_1}&\hdotsfor7&y_{k-1}&\SText{b_{k+1}}\\[1pt]
    \hline
    x_1&\hdotsfor1&x_{l-1}
    &\SText{a_{l+1}}&\hdotsfor1&a_{l^\prime}&\DText{x_{l^\prime}}&\hdotsfor1
        &x_{k-1}\\
    \DText{b_2}&\hdotsfor7&b_k&\SP
    \end{array}
  \end{equation*}
  Cutting at $(l,m)$ yields  
  \begin{equation*}\tag{$\tilde T$}
    \begin{array}{ccccccccccc}
    a_1&\hdotsfor2&a_l
    &\SText{a_{l+1}}&\hdotsfor1&a_{l^\prime}&\DText{x_{l^\prime}}&\hdotsfor1
        &x_{k-1}\\
    b_1&\DText{y_1}&\hdotsfor7&y_{k-1}
    &\SText{}\\[1pt]
    \hline
    x_1&\hdotsfor1&x_{l-1}
    &\SText{x_l}&\hdotsfor1&x_{l^\prime-1}&\DText{a_{l^\prime+1}}&\hdotsfor2
        &a_{k+1}\\
    \DText{b_2}&\hdotsfor7&b_k
    &\SText{b_{k+1}.}
    \end{array}
  \end{equation*}}

  \vbox{  
  \subsection*{$\mathbf{\text{Case }(11)\text{, }
      1=m^\prime\le l\le m\le l^\prime=k}$} The pair $T^\prime$ can be
  expressed in terms of the entries of $T$ as follows:
  \begin{equation*}\tag{$T^\prime$}
    \begin{array}{ccccccccccc}
    a_1&\hdotsfor2&a_l
    &\SText{x_l}&\hdotsfor4&x_{k-1}&\DText{a_{k+1}}\\
    b_1&\DText{y_1}&\hdotsfor4&y_{m-1}
    &\SText{b_{m+1}}&\hdotsfor2&b_{k+1}\\[1pt]
    \hline
    x_1&\hdotsfor1&x_{l-1}
    &\SText{a_{l+1}}&\hdotsfor4&a_k&\DText{}\\
    \DText{b_2}&\hdotsfor4&b_m
    &\SText{y_m}&\hdotsfor1&y_{k-1}.
    \end{array}
  \end{equation*}
  Cutting at $(l,m)$ yields  
  \begin{equation*}\tag{$\tilde T$}
    \begin{array}{ccccccccccc}
    a_1&\hdotsfor2&a_l
    &\SText{a_{l+1}}&\hdotsfor4&a_k&\DText{}\\
    b_1&\DText{y_1}&\hdotsfor4&y_{m-1}
    &\SText{y_m}&\hdotsfor1&y_{k-1}\\[1pt]
    \hline
    x_1&\hdotsfor1&x_{l-1}
    &\SText{x_l}&\hdotsfor4&x_{k-1}&\DText{a_{k+1}}\\
    \DText{b_2}&\hdotsfor4&b_m
    &\SText{b_{m+1}}&\hdotsfor2&b_{k+1}.
    \end{array}
  \end{equation*}}

  \vbox{
  \subsection*{$\mathbf{\text{B: }m\le l}$}
  In the Cases~$(7)$--$(9)$ and $(12)$ we have $m\le l$, thus the pair
  of two-rowed arrays $T=(T_1,T_2)\in\T^L_{k+1}\times\T^L_{k-1}$ looks
  like
  \begin{equation*}
    \begin{array}{cccccccccc}
      a_1&\hdotsfor5&a_l&\SText{a_{l+1}}&\hdotsfor1&a_{k+1}\\
      b_1&\hdotsfor2&b_m&\SText{b_{m+1}}&\hdotsfor4&b_{k+1}\\[1pt]
      \hline
      x_1&\hdotsfor4&x_{l-1}&\SText{x_l}&\hdotsfor1&x_{k-1}\\
      y_1&\dots&y_{m-1}&\SText{y_m}&\hdotsfor4&y_{k-1}.
    \end{array}
  \end{equation*}
  Cutting at $(l,m)$ we obtain $I(T)\in\T^L_k\times\T^L_k$:
  \begin{equation*}
    \begin{array}{cccccccccc}
      a_1&\hdotsfor5&a_l&\SText{x_l}&\dots&x_{k-1}\\
      b_1&\hdotsfor2&b_m&\SText{y_m}&\hdotsfor4&y_{k-1}\\[1pt]
      \hline
      x_1&\hdotsfor4&x_{l-1}&\SText{a_{l+1}}&\hdotsfor2&a_{k+1}\\
      y_1&\dots&y_{m-1}&\SText{b_{m+1}}&\hdotsfor5&b_{k+1}.
    \end{array}
  \end{equation*}}

  \vbox{
  \subsection*{$\mathbf{\text{Case }(7)\text{, }
      1=m\le l\le l^\prime\le m^\prime=k}$} The pair $T^\prime$ can be
  expressed in terms of the entries of $T$ as follows:
  \begin{equation*}\tag{$T^\prime$}
    \begin{array}{ccccccccccc}
    a_1&\hdotsfor2&a_l
    &\SText{x_l}&\hdotsfor1&x_{l^\prime-1}&\DText{a_{l^\prime+1}}&\hdotsfor2
        &a_{k+1}\\
    b_1
    &\SText{y_1}&\hdotsfor7&y_{k-1}&\DText{b_{k+1}}\\[1pt]
    \hline
    x_1&\hdotsfor1&x_{l-1}
    &\SText{a_{l+1}}&\hdotsfor1&a_{l^\prime}&\DText{x_{l^\prime}}&\hdotsfor1
        &x_{k-1}\\
    \SText{b_2}&\hdotsfor7&b_k&\DP
    \end{array}
  \end{equation*}
  Cutting at $(l,m)$ yields
  \begin{equation*}\tag{$\tilde T$}
    \begin{array}{ccccccccccc}
    a_1&\hdotsfor2&a_l
    &\SText{a_{l+1}}&\hdotsfor1&a_{l^\prime}&\DText{x_{l^\prime}}&\hdotsfor1
        &x_{k-1}\\
    b_1
    &\SText{b_2}&\hdotsfor7&b_k&\DText{}\\[1pt]
    \hline
    x_1&\hdotsfor1&x_{l-1}
    &\SText{x_l}&\hdotsfor1&x_{l^\prime-1}&\DText{a_{l^\prime+1}}&\hdotsfor2
        &a_{k+1}\\
    \SText{y_1}&\hdotsfor7&y_{k-1}&\DText{b_{k+1}.}
    \end{array}
  \end{equation*}}

  \vbox{
  \subsection*{$\mathbf{\text{Case }(8)\text{, }
      1=m\le l\le m^\prime\le l^\prime=k}$} The pair $T^\prime$ can be
  expressed in terms of the entries of $T$ as follows:
  \begin{equation*}\tag{$T^\prime$}
    \begin{array}{ccccccccccc}
    a_1&\hdotsfor2&a_l
    &\SText{x_l}&\hdotsfor4&x_{k-1}&\DText{a_{k+1}}\\
    b_1
    &\SText{y_1}&\hdotsfor4&y_{m^\prime-1}&\DText{b_{m^\prime+1}}&\hdotsfor2
        &b_{k+1}\\[1pt]
    \hline
    x_1&\hdotsfor1&x_{l-1}
    &\SText{a_{l+1}}&\hdotsfor4&a_k&\DText{}\\
    \SText{b_2}&\hdotsfor4&b_{m^\prime}&\DText{y_{m^\prime}}&\hdotsfor1
        &y_{k-1}.
    \end{array}
  \end{equation*}
  Cutting at $(l,m)$ yields  
  \begin{equation*}\tag{$\tilde T$}
    \begin{array}{ccccccccccc}
    a_1&\hdotsfor2&a_l
    &\SText{a_{l+1}}&\hdotsfor4&a_k&\DText{}\\
    b_1
    &\SText{b_2}&\hdotsfor4&b_{m^\prime}&\DText{y_{m^\prime}}&\hdotsfor1
        &y_{k-1}\\[1pt]
    \hline
    x_1&\hdotsfor1&x_{l-1}
    &\SText{x_l}&\hdotsfor4&x_{k-1}&\DText{a_{k+1}}\\
    \SText{y_1}&\hdotsfor4&y_{m^\prime-1}&\DText{b_{m^\prime+1}}&\hdotsfor2
        &b_{k+1}.
    \end{array}
  \end{equation*}}

  \vbox{
  \subsection*{$\mathbf{\text{Case }(9)\text{, }
      1=m\le m^\prime\le l\le l^\prime=k}$} The pair $T^\prime$ can be
  expressed in terms of the entries of $T$ as follows:
  \begin{equation*}\tag{$T^\prime$}
    \begin{array}{ccccccccccc}
    a_1&\hdotsfor5&a_l
    &\SText{x_l}&\hdotsfor1&x_{k-1}&\DText{a_{k+1}}\\
    b_1
    &\SText{y_1}&\hdotsfor1&y_{m^\prime-1}&\DText{b_{m^\prime+1}}&\hdotsfor5
        &b_{k+1}\\[1pt]
    \hline
    x_1&\hdotsfor4&x_{l-1}
    &\SText{a_{l+1}}&\hdotsfor1&a_k&\DText{}\\
    \SText{b_2}&\hdotsfor1&b_{m^\prime}&\DText{y_{m^\prime}}&\hdotsfor4
        &y_{k-1}.
  \end{array}
  \end{equation*}
  Cutting at $(l,m)$ yields
  \begin{equation*}\tag{$\tilde T$}
    \begin{array}{ccccccccccc}
    a_1&\hdotsfor5&a_l
    &\SText{a_{l+1}}&\hdotsfor1&a_k&\DText{}\\
    b_1
    &\SText{b_2}&\hdotsfor1&b_{m^\prime}&\DText{y_{m^\prime}}&\hdotsfor4
        &y_{k-1}\\[1pt]
    \hline
    x_1&\hdotsfor4&x_{l-1}
    &\SText{x_l}&\hdotsfor1&x_{k-1}&\DText{a_{k+1}}\\
    \SText{y_1}&\hdotsfor1&y_{m^\prime-1}&\DText{b_{m^\prime+1}}&\hdotsfor5
        &b_{k+1}.
  \end{array}
  \end{equation*}}

  \vbox{
  \subsection*{$\mathbf{\text{Case }(12)\text{, }
      1=m^\prime\le m\le l\le l^\prime=k}$} The pair $T^\prime$ can be
  expressed in terms of the entries of $T$ as follows:
  \begin{equation*}\tag{$T^\prime$}
    \begin{array}{ccccccccccc}
    a_1&\hdotsfor5&a_l
    &\SText{x_l}&\hdotsfor1&x_{k-1}&\DText{a_{k+1}}\\
    b_1&\DText{y_1}&\hdotsfor1&y_{m-1}
    &\SText{b_{m+1}}&\hdotsfor5&b_{k+1}\\[1pt]
    \hline
    x_1&\hdotsfor4&x_{l-1}
    &\SText{a_{l+1}}&\hdotsfor1&a_k&\DText{}\\
    \DText{b_2}&\hdotsfor1&b_m
    &\SText{y_m}&\hdotsfor4&y_{k-1}.
  \end{array}
  \end{equation*}
  Cutting at $(l,m)$ yields
  \begin{equation*}\tag{$\tilde T$}
    \begin{array}{ccccccccccc}
    a_1&\hdotsfor5&a_l
    &\SText{a_{l+1}}&\hdotsfor1&a_k&\DText{}\\
    b_1&\DText{y_1}&\hdotsfor1&y_{m-1}
    &\SText{y_m}&\hdotsfor4&y_{k-1}\\[1pt]
    \hline
    x_1&\hdotsfor4&x_{l-1}
    &\SText{x_l}&\hdotsfor1&x_{k-1}&\DText{a_{k+1}}\\
    \DText{b_2}&\hdotsfor1&b_m
    &\SText{b_{m+1}}&\hdotsfor5&b_{k+1}.
  \end{array}
  \end{equation*}}
\end{proof}
\newcommand{\cocoa} {\mbox{\rm C\kern-.13em o\kern-.07em C\kern-.13em
  o\kern-.15em A}}
\providecommand{\bysame}{\leavevmode\hbox to3em{\hrulefill}\thinspace}
\providecommand{\MR}{\relax\ifhmode\unskip\space\fi MR }
\providecommand{\MRhref}[2]{%
  \href{http://www.ams.org/mathscinet-getitem?mr=#1}{#2}
}
\providecommand{\href}[2]{#2}

\end{document}